\documentclass[a4paper,12pt]{article}
\usepackage {amsmath,amssymb}

\newcommand{\cl}[1]{\mathcal{#1}}

\newcommand{\F}{\mathbb{F}}

\newcommand{\Q}{\mathbb{Q}}

\newcommand{\beql}[1]{\begin{equation}\label{#1}}
\newcommand{\eeq}{\end{equation}}

\newcommand{\modd}[1]{\; ( \text{mod} \; #1)}
\newtheorem{theorem}{Theorem}
\newtheorem{defn}{Definition}

\newtheorem{corollary}{Corollary}
\newtheorem{lemma}{Lemma}

\begin{document}

\title{Irreducible polynomials over finite fields produced by
  composition of quadratics}

\author{D.R. Heath-Brown and G. Micheli\\Mathematical Institute, Oxford}

\date{}

\maketitle

\section{Introduction}

If $f(X)=f^{(1)}(X)=X^2+1$ then the iterates
$f^{(n)}(X)=f^{(n-1)}(f(X))$ are irreducible over $\Q$ for all
$n\ge 1$. Indeed the polynomials $f^{(n)}(X)$ are irreducible over $\F_3$.

This paper will explore a number of questions suggested by this
result.  In particular we investigate the situation in which one
composes more than one polynomial.  We make the following definition:

\begin{defn}
Let $f_1(X),\ldots,f_r(X)$ be polynomials of positive degree over a
finite field $\F_q$.  We say they are ``dynamically-irreducible'' if all
polynomials formed by composition of $f_1,\ldots,f_r$ are
irreducible over $\F_q$.
\end{defn}

Other authors have used the term ``stable'' for the situation $r=1$ (see also \cite{gomez2012stable}), but we believe
that ``dynamically-irreducible'' is more suggestive in this general context.

We will be particularly interested in the situation in which the $f_i$ are
quadratic polynomials. When $r=1$ a beautiful criterion for a  quadratic to be
dynamically-irreducible was developed by Boston and Jones
\cite[Proposition 2.3]{BJ}, but unfortunately the proof included a minor
error.\footnote{In the final display on page 1851, the second equality
  must be adjusted unless $g\circ f^{n-1}$ is monic.}  After
correcting this the criterion says that
$f(X)=a\{(X-b)^2+c\}\in\F_q$ is dynamically-irreducible if and only if $f$ is
irreducible, and $f^{(n)}(c)\not\in a\F_q^2$ for $n\ge 1$. For
example, when $f(X)=X^2+1\in\F_3$ we see that
$f^{(n)}(1)=2\not\in\F_3^2$, for all $n\ge 1$. This confirms the claim
made above that this polynomial is dynamically-irreducible (in fact, this is exactly the polynomial needed to build an infinite F-set in \cite{ferraguti2016existence} for $p=3$).

If $f(X)=(X-b)^2+b-2\in\F_q[X]$ we find that $f^{(n)}(b-2)=b+2$ for $n\ge 
1$. If $q$ is odd we therefore obtain a dynamically-irreducible polynomial
whenever $2-b$ and $2+b$ are both non-squares in $\F_q$. It is an
elementary exercise to show that a suitable $b$ may always be found.

Since the set of iterates above must eventually produce a cycle it is
clear from the above criterion that one can test in finite time
whether a given quadratic polynomial over $\F_q$ is
dynamically-irreducible.  Indeed Ostafe and Shparlinski \cite{OS} show that
one has a repetition within $O(q^{3/4})$ steps, so that one can test
whether a polynomial is dynamically-irreducible in $O(q^{3/4})$ operations.
\bigskip

The focus of this paper will be on large dynamically-irreducible sets of
quadratic polynomials. The criterion of Boston and Jones has been
extended as follows by Ferraguti, Micheli and Schnyder \cite[Theorem
  2.4]{FMS}.

\begin{lemma}\label{crit}
Let $f_i(X)=(X-b_i)^2+c_i\in\F_q$ be irreducible
polynomials for $1\le i\le r$.  Then the set $f_1,\ldots,f_r$ is
dynamically-irreducible if and only if every iterate
\beql{iter}
\left(f_{i_1}\circ f_{i_2}\circ \ldots \circ f_{i_n}\right)(c_j)
\eeq
with $n\ge 1$ and $1\le i_1,\ldots,i_n,j\le r$, is a non-square
in $\F_q$.
\end{lemma}

Notice that this is only formulated for the case in which all the
polynomials $f_i$ are monic. However, it has been pointed out by Alina Ostafe that one should be able to extend this criterion to non-monic quadratics, with obvious consequences for the various results in the present paper. Indeed she is able to exhibit 
large dynamically-irreducible sets of non-monic quadratics over fields of prime 
cardinality $p\equiv 1\modd{4}$. 
As an example (taken from \cite{FMS}) we
may consider the polynomials $(X-a)^2+a$ and $(X-a-1)^2+a$ in $\F_q$,
with $q\equiv 1\modd{4}$.  One then finds that the iterates
(\ref{iter}) take only the values $a$ and $a+1$.  Thus, choosing $a$
so that both $a$ and $a+1$ are non-squares (as we always can) we
obtain a dynamically-irreducible set of size 2.

Our first new results concern the existence of large dynamically-irreducible
sets of quadratic polynomials.  Of course, without the restriction to
quadratic polynomials, one can have infinitely large sets over $\F_q$,
since once $f$ is dynamically-irreducible, the set $f^{(1)},f^{(2)},\ldots$
is dynamically-irreducible. We write $M(q)$ for the size of the maximal set
of monic quadratic polynomials over $\F_q$ which is dynamically-irreducible.

\begin{theorem}\label{constr}
Let $p$ be a prime and $h$ a non-zero element of $\F_p$. Then the polynomial 
$X^p-X-h$ is irreducible over $\F_p$. Take $\xi$ to be a root of  
$X^p-X-h$ let $K=\F_p(\xi)$, so that $K$ is a finite field with
$q=p^p$ elements. Define polynomials
\[f_{b,c}(X)=(X-b-\xi)^2+c+\xi\in K[X],\;\;\mbox{for
  all}\;(b,c)\in\F_p^2.\]
If $p\equiv 1\modd{4}$,  then the polynomials $f_{b,c}(X)$ form a dynamically-irreducible set over $K$ whenever
$h$ is a non-square in $\F_p$ .
\end{theorem}

\begin{theorem}\label{count1}
  There are infinitely many finite fields for which
\[M(q)\ge\tfrac12(\log q)^2.\]
\end{theorem}

 Notice that the explicit construction in Theorem \ref{constr} 
yields $M(p^p)\ge p^2$ and hence provides examples with
$M(q)\ge (\log q)^2(\log\log q)^{-2}$.

Unfortunately Theorems \ref{constr} and \ref{count1} do not provide
fields of prime order in which there are large dynamically-irreducible sets.
A little experimentation shows that $M(3)=1$, $M(5)=3$, $M(7)=2$, and
$M(11)=1$.  In particular, over $\F_5$, the three polynomials
\[(x-2)^2+2,\; (x-3)^2+2,\; x^2+3\]
form a dynamically-irreducible set. Similarly, over $\F_{13}$, the three polynomials 
\[(x - 1)^2 - 2,\; (x - 9)^2 -6,\; (x - 3)^2 - 5\]
form a dynamically-irreducible set.  As mentioned above, we can find a
dynamically-irreducible set of size 2 over $\F_p$ whenever $p\equiv 1\modd{4}$.

Our next result describes the set of iterates (\ref{iter}).
\begin{theorem}\label{ci}
Let $f_1,\ldots,f_r$ be a dynamically-irreducible set of monic quadratic
polynomials over a finite field $\F_q$ of odd characteristic, and
suppose that $r\ge 2$.  Then
the set of iterates (\ref{iter}), together with the values of the $c_j$,
has size at most
$4(\log q)^2\sqrt{q}$, uniformly in $r$.
\end{theorem}

One should compare this to the corresponding result for $r=1$ given by
Ostafe and Shparlinski \cite{OS} which we mentioned above, and in
which the size is $O(q^{3/4})$. It may seem counter-intuitive that one
should have fewer values (\ref{iter}), despite having more polynomials
to use. However the requirement that all the elements obtained should
be non-squares has the over-riding effect.

We have three corollaries to Theorem \ref{ci}.

\begin{corollary}\label{count2}
For any odd prime power $q$ we have $M(q)\le 32(\log q)^4q$.
\end{corollary}

\begin{corollary}\label{corollaryDS1}
Let $f_1,\dots,f_r$ be a set of monic dynamically-irreducible polynomials, all with the same value for $c_j$.
Then $r\leq 8 q^{1/2}(\log q )^2$.
\end{corollary}

\begin{corollary}\label{alg}
  Let $\F_q$ be a finite field of odd characteristic.
There is an algorithm to test whether or not a set
$f_1,\ldots,f_r$ of monic quadratic
polynomials over $\F_q$ is dynamically-irreducible, which takes $O(r(\log
q)^3q^{1/2})$ operations, and requires $O((\log q)^2q^{1/2})$ storage
locations. 
\end{corollary}

The upper bound in Corollary \ref{count2} is disappointingly weak.
The proof we give is extremely simple and discards much of the
available information.  It would be interesting to know whether a more
sophisticated approach would lead to an improved estimate for $M(q)$.

As mentioned above
Ostafe and Shparlinski \cite{OS} show that one can test
whether a quadratic polynomial is dynamically-irreducible in $O(q^{3/4})$ operations.
It is pleasing to see that Corollary \ref{alg} produces a faster
algorithm for $r=2$, say, than one has for $r=1$. However when $r=1$
the space requirement is $O(1)$, while our algorithm needs non-trivial
amounts of memory since it uses a tree structure.
\bigskip

{\bf Acknowledgments.} Giacomo Micheli was supported by the Swiss National Science Foundation grant number 161757.

We would also like to record our thanks to Professor Edith Elkind, who
suggested to us the use of Red-Black trees in the proof of Corollary \ref{alg}.

\section{Proof of Theorems \ref{constr} and \ref{count1}}

To prove that $T(X)=X^p-X-h$ is irreducible over $\F_p$
we begin by observing that if $\xi$ is a root of $T(X)$ then so is
$\xi+a$ for any $a\in\F_p$.  Thus the roots of $T(X)$ are precisely
the values $\xi+a$ as $a$ runs over $\F_p$.  Let $G(X)\in\F_p[X]$
be a factor of
$T(X)$ of degree $d$, with roots $\xi+a_1,\ldots, \xi+a_d$ say. Then
the sum of the roots is given by the coefficient of $X^{d-1}$, and
hence lies in $\F_p$.  However the sum of the roots will be
$d\xi+a_1+\ldots+a_d$. If this is in $\F_p$ then either $p\mid d$ or
$\xi\in\F_p$. However $T(n)=n^p-n-h=-h\not=0$ for every $n\in\F_p$,
whence $G$ must have degree $0$ or $p$.

We now come to the key idea for the two theorems. Polynomials over a
field form a (non-commutative) semigroup under composition.  The
polynomial $X$ acts as identity element, and the polynomial $X+j$,
where $j$ is constant, has a two sided inverse $X-j$.  Suppose now
that we have an extension $M/L$ of finite fields of odd
characteristic. Let $g_{b,c}(X)=(X-b)^2+c$, where $b,c$ run over $L$,
and let $\ell(X)=X+\xi$ for some $\xi\in M$. Then, given any composition 
\[G=g_{b_1,c_1}\circ \ldots\circ g_{b_n,b_n}\]
we will have
\[F=f_{b_1,c_1}\circ \ldots\circ f_{b_n,b_n}\]
with 
\[F(X)=(\ell\circ G\circ\ell^{-1})(X)=G(X-\xi)+\xi\]
and
\[f_{b,c}(X)=(\ell\circ g_{b,c}\circ\ell^{-1})(X)=(X-b-\xi)^2+c+\xi.\]
Here $f_{b,c}$ is irreducible over $M$ if $-c-\xi$ is not a square in $M$.  
For any $c\in L$ we have $F(c+\xi)=G(c)+\xi$.  However
$G(c)\in L$. It therefore follows from Lemma \ref{crit} that the
polynomials $f_{b,c}$ form a dynamically-irreducible set over $M$, 
provided that every
element of the additive cosets $L+\xi$ and $L-\xi$ is a non-square in $M$. 
If $-1$ is a square in $L$ it will be enough to consider $L+\xi$.

For Theorem \ref{constr} we take $L=\F_p$ and $M=K$.  If $a+\xi$ were a square in 
$K$, then its norm $N_{K/\F_p}(a+\xi)$ would be a square in $\F_p$.
However this norm is simply the product of the roots of
$F(X)=X^p-X-h$, which is $h$.  Thus if $h$ is a non-square in $\F_p$,
then every element $a+\xi$ is a non-square in $K$.  This completes the
proof of Theorem \ref{constr}.
\bigskip

For Theorem \ref{count1}, we let $q$ be a power $p^e$ of a prime
$p\equiv 1\modd{4}$,
and apply the previous ideas with $L=\F_p$ and $M=\F_q$.  We
therefore hope to find an element $\alpha\in\F_q$ such that $a+\alpha$
is a non-square in $\F_q$ for every $a\in\F_p$.  For this we will
use estimates for character sums.

Let $\chi$ be the quadratic character for $\F_q$, and consider
\beql{Sd}
S=\sum_{\alpha\in\F_q}\prod_{a\in\F_p}\{1-\chi(a+\alpha)\}.
\eeq
If $\alpha\in\F_p$ the factor corresponding to $a=1-\alpha$ is zero,
so that such an $\alpha$ makes no contribution to $S$. Thus if $S>0$ there 
must be some $\alpha\in \F_q\setminus \F_p$ such that $a+\alpha$ is a non-square for every
$a$. Our goal is therefore to show that $S>0$.

On expanding the product in (\ref{Sd}) we see that
\[S=q+\sum_{H}\varepsilon(H)\sum_{\alpha\in\F_q}\chi(H(\alpha)),\]
where $H$ runs over polynomials of the form
\[H(X)=\prod_{a\in A}(X+a)\]
for the various non-empty subsets $A\subseteq \F_p$, and 
$\varepsilon(H)=(-1)^{\# A}$.  It follows from Weil's bound for
character sums that
\[\left|\sum_{\alpha\in\F_q}\chi(H(\alpha))\right|\le {\rm
  deg}(H)q^{1/2}\le q^{1/2}p,\]
whence
\[S\ge q-q^{1/2}p(2^p-1).\]
It follows that a suitable $\alpha$ exists as soon as $q\ge p^24^p$. In
particular, it suffices to take $e$ as the first integer greater than equal
to $2+p(\log 4)/(\log p)$. With this choice we have
\[e\le 3+p(\log 4)/(\log p)\le \sqrt{2}\frac{p}{\log p}\]
provided that we take $p$ large enough. As a result we have $\log q\le
\sqrt{2} p$, whence $M(q)\ge p^2\ge \tfrac{1}{2}(\log q)^2$, as claimed.

\section{Proof of Theorem \ref{ci}}

Let $\mathcal{J}\subseteq\F_q$ denote the set of values
(\ref{iter}), together with the values of the $c_j$.
Assuming that $f_1,\ldots,f_r$ form a dynamically-irreducible
set, one sees that
$\left(f_{i_1}\circ f_{i_2}\circ \ldots \circ f_{i_n}\right)(j)$
must be a non-square in $\F_q$ for every $j\in\cl{J}$. We shall work
with a finite subset $\cl{F}$ of distinct polynomials given as compositions
$f_{i_1}\circ\ldots\circ f_{i_n}$ for some fixed $n\ge 1$. Thus 
$\cl{F}$ consists of distinct irreducible monic polynomials of degree $2^n$.

If $\chi$ is the quadratic
character for $\F_q$ we have
\[\prod_{F\in\cl{F}}\{1-\chi(F(j))\}\left\{\begin{array}{cc}
=2^{\#\cl{F}}, & j\in\cl{J},\\ \ge 0, &
j\not\in\cl{J},\end{array}\right.\]
and it follows that
\beql{bd}
\#\cl{J}\le 2^{-\#\cl{F}}\sum_{j\in\F_q}\prod_{F\in\cl{F}}\{1-\chi(F(j))\}.
\eeq
We proceed to expand the product and to use the Weil bound to
estimate the resulting character sums
\[\sum_{j\in\F_q}\chi(F_1(j)\ldots F_m(j)).\]
The polynomial $F_1\ldots F_m$ will be square-free, with degree
$2^n m\le 2^n\#\cl{F}$, whence the Weil bound produces
\[\left|\sum_{j\in\F_q}\chi(F_1(j)\ldots F_m(j))\right|\le
2^n(\#\cl{F})q^{1/2}\] 
as long as $m\ge 1$. The sum corresponding to $m=0$ is just $q$, and
there are $2^{\#\cl{F}}-1$ other sums, so that (\ref{bd}) yields
\beql{wl}
\#\cl{J}\le 2^{-\#\cl{F}}q+2^n(\#\cl{F})q^{1/2}.
\eeq

However at this point we encounter a potential difficulty.  If it were
true that all compositions $f_{i_1}\circ\ldots\circ f_{i_n}$ were
different we could take $\#\cl{F}=2^n$. Unfortunately this is not the
case. For example if $f_1(X)=(X+1)^2, f_2(X)=X^2$
and $f_3(X)=X^2+1$ then $f_1\circ f_2=f_2\circ f_3$. Of course this
problem does not arise when $r=1$ since in this situation the iterates will have
different degrees. We will show that the case $r=2$ is also
satisfactory.

\begin{lemma}\label{free}
Let $f_1,f_2$ be distinct monic quadratic polynomials over $\F_q$, with $q$ odd.
Suppose that
\beql{eq}
f_{i_1}\circ\ldots\circ f_{i_n}=f_{j_1}\circ\ldots\circ f_{j_m}
\eeq
with $i_1,\ldots,i_n,j_1,\ldots,j_m\in\{1,2\}$.  Then $m=n$ and
$i_h=j_h$ for every index $h$.
\end{lemma}

We will prove this in a moment, but first we use it to complete the
proof of Theorem \ref{ci}. We 
take $\cl{F}$ to consist of all compositions of $n$ polynomials each
of which is either $f_1$ or $f_2$.  By Lemma \ref{free} we obtain
$2^n$ distinct polynomials this way, and (\ref{wl}) becomes
\[\#\cl{J}\le 2^{-2^n}q+4^nq^{1/2}.\]
Finally, we choose $n\ge 1$ so that
\[\frac{\log q}{\log 4}\le 2^n <2\frac{\log q}{\log 4},\]
whence
\[\#\cl{J}\le \sqrt{q}+\left(\frac{\log q}{\log 2}\right)^2\sqrt{q}\le
4(\log q)^2\sqrt{q},\]
as claimed.
\bigskip

It remains to establish Lemma \ref{free}. If (\ref{eq}) holds, the
two sides have degrees $2^n$ and $2^m$ so that we must have $n=m$. We
now argue by contradiction, supposing that we have a non-trivial
relation (\ref{eq}) in which $n$ is minimal. Then
$f_i\circ F=f_j\circ G$, say, in which either $f_i\not=f_j$, or
$f_i=f_j$ but $F\not=G$.  Let $f_i(X)=(X-a)^2+b$ and
$f_j(X)=(X-c)^2+d$. Then
\begin{eqnarray*}
b-d&=&\left(G(X)-c\right)^2-\left(F(X)-a\right)^2\\
&=&\{G(X)+F(X)-a-c\}\{G(X)-F(X)+a-c\}.
\end{eqnarray*}
Since $F$ and $G$ are monic, and $\F_q$ has odd characteristic,
the polynomial $G(X)+F(X)-a-c$ has positive degree.  
We therefore see
that $b=d$, and that $G(X)-F(X)+a-c=0$. If $a=c$ we would
have $f_i=f_j$ and $F=G$, giving us a contradiction.  Hence $a\not=c$
so that $f_1$ and $f_2$ are $(X-a)^2+b$ and $(X-c)^2+b$, in some order.
Moreover $F\not=G$, so that we must have $n\ge 2$.

Now let $F=f_r\circ U$ and $G=f_s\circ V$, say, with
$f_r(X)=(X-e)^2+b$ and $f_s(X)=(X-f)^2+b$ for appropriate values
$e,f=a$ or $c$.  Then
\begin{eqnarray*}
0&\not=&a-c\\
&=&F(X)-G(X)\\
&=&\left(U(X)-e\right)^2-\left(V(X)-f\right)^2\\
&=&\{U(X)+V(X)-e-f\}\{U(X)-V(X)+f-e\}.
\end{eqnarray*}
This however is impossible, since the factor $U(X)+V(X)-e-f$ has
positive degree, its leading coefficient being $2\not=0$. This
contradiction completes the proof of the lemma.

\section{Proof of the Corollaries}

For the proof of Corollary \ref{count2} it will be notationally
convenient to write $B=[4(\log q)^2\sqrt{q}]$.
We begin by observing that the
coefficients $c_1,\ldots,c_r$ contain at most $B$
distinct values, by Theorem \ref{ci}.  Moreover as $i$ runs from 1 to
$r$ we obtain at most $B$ distinct values for
$f_i(c_1)=(c_1-b_i)^2+c_i$. We now sort the polynomials $f_i$
according to the value taken by $c_i$.  For polynomials $f_i$ with
a given value $c_i=c$ there are at most $B$ values that
$f_i(c_1)=(c_1-b_i)^2+c$ can take, each of which corresponds to at most
$2B$ possible choices for $b_i$.  Thus there are at most $2B$
polynomials with any given value for $c_i$, and at most $2B^2$
polynomials in total.  This proves the corollary.
\bigskip

We now prove Corollary \ref{corollaryDS1}. By the assumption, we have that $c_j=u\in \F_q$ for all $j\in \{1,\dots, r\}$. Consider the map 
\[\Gamma: \{f_1,\dots f_r\}\longrightarrow \F_q\]
\[f_i\mapsto f_i(u).\]
Suppose that $f_i(u)=f_j(u)=f_k(u)$ for distinct $i,j,k\in \{1,\dots,r\}$, then
\[(u-a)^2=(u-b)^2=(u-c)^2\]
for some $a,b,c\in \mathbb F_q$. This implies $u=(a+b)/2=(a+c)/2$, which forces $b=c$, but then $f_j=f_k$. This shows that the map $\Gamma$ is at most $2$-to-1.
Since the image of $\Gamma$ is clearly contained in the set of iterates (\ref{iter}),  the claim follows by applying directly Theorem \ref{ci}.

\bigskip

Finally we tackle Corollary \ref{alg}.  We assume that the elements of
$\F_q$ are described in such a way as to enable us to impose an
ordering on them.  For example, if we have $q=p^d$ and take
$\F_q=\F_p[\theta]$ where $\theta$ is the root of an irreducible
polynomial of degree $d$, we might arrange the elements
$a_0+a_1\theta+\ldots+a_{d-1}\theta^{d-1}$ using lexicographic order.
Our algorithm inputs a set of irreducible quadratic polynomials
$f_1,\ldots,f_r$.  It then uses the following steps.

\begin{enumerate}
\item[1.] Build a Red-Black tree $T$, starting with the empty set and
  successively adding the values $c_1,\ldots,c_r$, using the ordering
  mentioned above, and discarding any duplicate values.
\item[2.] Set $j=0$ and form a list $L_j$ of the elements in $T$.
\item[3.] Set $L=\emptyset$.  Then, for each $i=1,\ldots,r$ and each
  $t\in L_j$
  \begin{enumerate}
  \item[i)] Check whether $f_i(t)$ is a square.  If it is, STOP
    (the polynomials do not form a dynamically-irreducible set).
\item[ii)] Otherwise check whether $f_i(t)$ is already in $T$.  If it
  is, move on to the next pair $i,t$.
  \item[iii)] If $f_i(t)$ is not yet in $T$ insert it into the tree,
    and add it to the list $L$.
  \end{enumerate}
  \item[4.] If $L$ is empty, STOP (the polynomials form a
  dynamically-irreducible set).
\item[5.] If $L$ is non-empty, add one to the value of $j$ and set
  $L_j=L$.  Return to Step 3.
\end{enumerate}

If the algorithm stops at Step 3(i) it has found a square iterate, and
so it correctly reports that we do not have a dynamically-irreducible
set. Otherwise, when we begin Step 3 the tree $T$ contains all
distinct values (\ref{iter}) with $n\le j$, and the list $L_j$ will
consist of all such values with $n=j$ which cannot be obtained from
any smaller $n$. If we stop at Step 4 then there are no new values
with $n=j+1$, so that further iteration will always produce results
already contained in $T$. In this case all iterates will be
non-squares and the algorithm correctly reports that we have a
dynamically-irreducible set of quadratics.

One readily sees that, whenever we begin Step 3, the size of $T$ is
$\# L_0+\ldots+\# L_j$.  On the other hand, Theorem \ref{ci} shows
that $T$ can have at most $B=[4(\log q)^2\sqrt{q}]$ elements. Since
the algorithm will terminate at Step 4 unless $L$ is non-empty, it is
clear that we must stop after at most $B$ loops.

To analyze the running time of the algorithm we note that we may test
an element of $\F_q$ to check whether it is a square using $O(\log q)$ field
operations. Moreover, since $T$ has size at most $B$, we can check
whether an element belongs to $T$, and if not add it to $T$, in
$O(\log q)$ operations. Step 3 requires $r(\# L_j)$ tests of this
type, so that the total number of operations needed is
\[\ll r(\# L_0+\# L_1+\ldots)\log q\ll rB\log q\ll r(\log
q)^3q^{1/2},\]
as claimed. As to the memory requirement, the tree $T$ and the lists
$L_j$ will need space $O(B)$.  This completes the proof of the corollary.

The use of a Red-Black tree, into which one may insert new elements in
order, was suggested to us by Professor Edith Elkind.  It is a
pleasure to record our thanks for this.

\bigskip

\bigskip

Mathematical Institute,

Radcliffe Observatory Quarter,

Woodstock Road,

Oxford

OX2 6GG

UK

\bigskip

{\tt rhb@maths.ox.ac.uk}\;\; and\;\; {\tt giacomo.micheli@maths.ox.ac.uk}

\end{document}